\documentclass[12pt]{amsart}
\usepackage{amssymb}
\usepackage{color}

\allowdisplaybreaks

\makeatletter
\@namedef{subjclassname@2020}{%
	\textup{2020} Mathematics Subject Classification}
\makeatother

\newtheorem{theorem}{Theorem}
\newtheorem{lemma}{Lemma}
\newtheorem{remark}{Remark}
\newtheorem{corollary}{Corollary}

\begin{document}
	
\title[Matrix-transform means\dots]{Pointwise Convergence and Approximation of Matrix-Transform Means of Vilenkin--Fourier Series}

\author[N. Anakidze]{N. Anakidze}
\address {Nino Anakidze, The University of Georgia, School of Science and Technology, 77a Merab Kostava St, Tbilisi, 0171, Georgia.}
\email{nino.anakidze@mail.ru}

\author[M. Szilveszter]{M. Szilveszter}
\address{M\'at\'e Szilveszter, Department of Numerical Analysis, E\"otv\"os Lor\'and University, Budapest, Hungary, P\'azm\'any P\'eter s\'et\'any 1/C, H-1117.}
\email{szilveszter.mate@inf.elte.hu}

\author[G. Tephnadze]{G. Tephnadze}
\address {George Tephnadze, The University of Georgia, School of Science and Technology, 77a Merab Kostava St, Tbilisi, 0171, Georgia.}
\email{g.tephnadze@ug.edu.ge}

\author[F. Weisz]{F. Weisz}
\address{Ferenc Weisz, Department of Numerical Analysis, E\"otv\"os Lor\'and University, Budapest, Hungary,
P\'azm\'any P\'eter s\'et\'any 1/C, H-1117.}
\email{weisz@inf.elte.hu}

\date{}

\begin{abstract} 
We study convergence, divergence, and approximation of matrix-transform means of Vilenkin--Fourier series. For bounded groups, we prove convergence along ${M_n}$ at Vilenkin--Lebesgue points without additional diagonal conditions. For arbitrary bounded or unbounded groups, every regular matrix-transform method admits divergence on any prescribed null set along any increasing sequence. We also obtain $L_p(G_m)$ approximation estimates in terms of the modulus of continuity.
\end{abstract}

\subjclass[2020]{42C10, 40C05, 42A20, 41A25}
\keywords{Fourier series, Vilenkin system, rate of approximation, matrix-transform means}

\maketitle

\section{Notation and definitions}

Let $\mathbb{P}$ be the set of positive natural numbers and $\mathbb{N}:=\mathbb{P}\cup \{0\}$. Let $m:=(m_{0}, m_{1},\cdots)$ denote a sequence of integers, each of which is at least $2$. Define the group $G_{m}$ as the complete direct product of the additive group $
\mathbb{Z}_{m_{k}}:=\{0,1,\cdots,m_{k}-1\}
$ of integers modulo $m_{k}$. The group operation is addition modulo $m_k$. 
Equip each $\mathbb{Z}_{m_k}$ with the discrete topology.
The normalized Haar measure $\mu$ on $\mathbb{Z}_{m_k}$ is given by 
$\mu (\{ 0\})=\mu (\{ 1\})=\cdots=\mu (\{ m_k-1\})=1/m_k,$
i.e., the measure of a singleton is $1/m_k$. 
The group $G_m:=
\overset{\infty}{\underset{k=0}{{\times}}} \mathbb{Z}_{m_k}$
is called the Vilenkin group. 

In this paper, we consider both bounded and unbounded Vilenkin groups. Whenever we restrict our attention to bounded Vilenkin systems, we explicitly assume that
$2\leq m_k\leq R<\infty,$ $k\in\mathbb N,$
where
$R:=\sup_{k\in\mathbb N}m_k.$

The elements of $G_{m}$ are represented by sequences
$
x:=\left( x_{0},x_{1},\cdots,x_{j},\cdots\right), \ \left( x_{j}\in \mathbb{Z}_{m_{j}}\right) .
$
Let
$
e_j:=\left(0,0,\cdots,1,0,\cdots\right),
$
where the $j$-th coordinate is equal to $1$ and all other coordinates are equal to $0.$
The group operation on $G_m$ is the coordinate-wise addition (denoted by $+$),
the normalized Haar measure  $\mu $ is the product measure and the topology is the product topology (see \cite{PTW}).

For $x\in G_m$ and $n\in\mathbb{P}$, define
\[
I_0(x):=G_m,\qquad
I_n(x):=\{y\in G_m :
y=(x_0,\dots,x_{n-1},y_n,y_{n+1},\dots)\}.
\]
Thus, the coordinates $0,\dots,n-1$ are fixed, while the coordinates with indices
$k\geq n$ are arbitrary. These intervals form a neighborhood base in $G_m$. 

Let $I_{n}:=I_{n}(0),$ that is, the first $n$ coordinates of every $y\in I_n$ are equal to zero, while the remaining coordinates are arbitrary.

Let $L_p(G_m)$, $1\leq p<\infty$, denote the usual Lebesgue spaces on $G_m$, equipped with the norm $\|\cdot\|_p$. For $p=\infty$, we work with the space $C(G_m)$ of continuous functions on $G_m$, endowed with the supremum norm
$
\|f\|_\infty:=\sup_{x\in G_m}|f(x)|.
$
Thus, whenever $p=\infty$ occurs below, $L_\infty(G_m)$ is understood to mean $C(G_m)$ equipped with the supremum norm.

We now recall some basic notions of Vilenkin-Fourier analysis.
If we define the so-called generalized number system based on $m$ by
$
M_{0}:=1,\ M_{k+1}:=m_{k}M_{k},$ $(k\in \mathbb{N}),$
then every natural number $n$ can uniquely be expressed in the form $n=\sum_{j=0}^{\infty }n_{j}M_{j},$ where $n_{j}\in \mathbb{Z}_{m_{j}}$ $(j\in\mathbb{N})$ and only finitely many digits $n_{j}$ are nonzero. 
Let the order of $n\in\mathbb{P}$ be denoted by $\vert n\vert :=\max \{ j\in\mathbb{N}: n_j\neq 0\}.$  It means $M_{|n|}\leq n< M_{|n|+1}$. 

The Rademacher functions are defined as
\begin{equation*}
r_{k}\left( x\right) :=\exp \left( 2\pi \imath x_{k}/m_{k}\right) ,\text{ }\left(\imath^{2}=-1, \ x\in G_{m},\ k\in\mathbb{N}\right) .
\end{equation*}
The Vilenkin system is the product system generated by the Rademacher functions
\begin{equation*}
\psi _{n}(x):=\prod\limits_{k=0}^{\infty }r_{k}^{n_{k}}\left( x\right) \,\,\ \ \,\left( n\in \mathbb{N}\right).
\end{equation*}

If $m_k=2$ for every $k\in\mathbb{N}$, then $G_m$ is the dyadic group, and the Vilenkin system coincides with the Walsh system.

Let ${\mathcal P}_n$ be the collection of Vilenkin polynomials of order less than $n$, that is, functions of the form
$
P(x)=\sum_{k=0}^{n-1}a_k\psi_k(x),
$
where $n\in\mathbb{P}$ and $\{a_k\}$ is a sequence of complex numbers.
The collection of all Vilenkin polynomials is denoted by
$
{\mathcal P}:=\bigcup_{n\in\mathbb{P}}{\mathcal P}_n.
$

It is known \cite{hew} that the Vilenkin system $(\psi_n, n\in\mathbb{N})$ is the character system of $(G_m,+)$.
For any $f\in L_{p}(G_m)$ and  $\delta>0$, the  modulus of continuity is defined by
\[
\omega_p(f,\delta):=\sup_{|t|<\delta}\left\Vert f(\cdot-t)-f(\cdot)\right\Vert_p,
\]
 where  $-$ is the inverse operation of  $+$ defined on $G_m$ and 
$\lvert x\rvert:=\sum_{i=0}^\infty {x_i}/{M_{i+1}}$ for all $x\in G_m.$
In the case $f\in C(G_m)$ we replace $p$ by $\infty$.

Fourier coefficients, partial sums, and Dirichlet kernels are
defined by
\begin{equation*}
	\hat f(j):=\int_{G_m} f(x)\overline{\psi}_{j}(x)d\mu(x),
\ \ \ \
S_{k}(f):=\sum_{j=0}^{k-1}\hat f(j)\psi_{j}, \ \ \ \
D_{n}:=
\sum_{k=0}^{n-1}\psi_{k},\ D_{0}:=0.
\end{equation*}
Let $n\in\mathbb{N}$. Then we have (see  \cite{AVD} and \cite{AG})
\begin{equation}\label{dn40}
\int_{G_m}\left\vert D_{M_n}(x)\right\vert d\mu(x)=\int_{G_m} D_{M_n}(x) d\mu(x)=1, \quad
\quad \hspace*{0in}D_{M_{n}}\left( x\right) =\left\{
\begin{array}{l}
\text{ }M_n,\text{\thinspace \thinspace \thinspace  if\thinspace
	\thinspace }x\in I_{n}, \\
\text{ }0,\text{\thinspace \thinspace \thinspace \thinspace \thinspace \thinspace \thinspace \thinspace if
	\thinspace \thinspace }x\notin I_{n},
\end{array}
\right.
\end{equation}

\begin{eqnarray}\label{Dbont}
D_{M_n-i}
&=&D_{M_n}-\psi_{M_n-1}\overline{D_{i}},\,\,\ 0\le i< M_n.
\end{eqnarray}

Fej\'er means and Fej\'er kernels are defined as
\[
\sigma_{n}(f):=\frac{1}{n}\sum_{k=1}^{n}S_{k}(f)
\ \ \ 
\text{and} \ \ \ 
K_{n}:=\frac{1}{n}\sum_{k=1}^{n}D_{k}.
\]

Let $n\in\mathbb P$. Then  (see \cite{G1} and \cite{PTW1})
\begin{eqnarray} \label{fn40}
\int_{G_m} K_n (x)d\mu(x)=1,  \ \ \ \ \ 
\sup_{n\in\mathbb{P}}\int_{G_m}\left\vert K_n(x)\right\vert d\mu(x)\leq 4R^3.
\end{eqnarray}
and
\begin{eqnarray} \label{fn41}
n|K_n|
\leq
c\sum_{A=0}^{|n|}M_A|K_{M_A}| \ \ \text{and} \ \
|K_{M_n}(u)|
\leq
c\sum_{s=0}^{n}M_s
\sum_{r=1}^{m_s-1}
1_{I_n(u-r e_s)},
\end{eqnarray}

Let $\{ q_k:k\in\mathbb{N}\}$ be a sequence of non-negative numbers. The $n$-th Nörlund mean and the $n$-th $T$-mean are defined by
\begin{equation*} 
t_{n}(f):=\frac{1}{Q_n}\sum_{k=1}^{n}q_{n-k}S_k(f)
\ \ \ \text{	and} \ \ \ 
T_n(f):=\frac{1}{Q_n}\overset{n}{\underset{k=1}{\sum }}q_{k-1}S_{k}(f), 
\end{equation*}
where 
$Q_n:=\sum_{k=0}^{n-1}q_k, \ (n\in\mathbb{P}).$
It is  always assumed that $q_0>0$ and 
\begin{equation}\label{Nor}
	\lim_{n\to\infty}Q_n=\infty.
\end{equation}
In this case, the N\"orlund method generated by  $\{ q_k:k\in\mathbb{N}\}$ is regular if and only if
\begin{equation}\label{reg}
	\lim_{n\to\infty}\frac{q_{n-1}}{Q_{n}}=0.
\end{equation}
For these notions and results, we refer the reader to \cite{Moore}, pp. 37-38.

We also note that the $T$-summability method   generated by  $\{ q_k:k\in\mathbb{N}\}$ is regular if and only if condition \eqref{Nor} is fulfilled (for details, see \cite{PTW}).


Let
$
T:=\left(t_{i,j}\right)_{i,j=1}^{\infty}
$
be an infinite triangular matrix. The $n$-th matrix-transform mean associated with $T$ is defined by
\[
\sigma_{n}^{T}(f):=\sum_{i=1}^{n}t_{i,n}S_{i}(f),
\]
where $\{ t_{i,n}: 1\leq i\leq n, \ i\in\mathbb{P}\}$ is a finite sequence of non-negative numbers and
\begin{equation}\label{28}
\sum_{i=1}^n t_{i,n}=1 \qquad (n\in\mathbb{P}).
\end{equation}

By the Silverman--Toeplitz theorem (see \cite[pp.~43--48]{Hardy}),
the matrix-transform summability method is regular if and only if
\begin{equation}\label{29}
\lim_{n \to \infty} t_{k,n} =0, \quad \text{for each fixed } k.
\end{equation}

The matrix-transform summation method is a common generalization of several well-known summation methods. It is immediate that the N\"orlund means, the $T$-means, the  Fej\'er means, the $(C,\alpha)$ means and the N\"orlund and Riesz logarithmic  means are special cases of the matrix-transform summation method. For example, if we suppose \eqref{28} and
\begin{equation*}
t_{k,n}:=\frac{q_{n-k}}{Q_{n}}, \ \ \ \ 	t_{k,n}:=\frac{q_{k-1}}{Q_{n}} \ \ \ \ (k=1,2,\cdots,n),
\end{equation*}
then we get the N\"orlund and $T$-means, respectively.
It is easily seen that \[ \sigma_n^{T}(f;x) =\int_{G_m} f(x-u)K_n^{T}(u)\,d\mu(u), \qquad x\in G_m, \] where $ K_n^{T}:=\sum_{k=1}^{n}t_{k,n}D_k $ is the kernel associated with the matrix-transform mean $\sigma_n^{T}$. Thus, the study of matrix-transform means can be reduced to obtaining appropriate estimates for the kernels $K_n^{T}$.
Let 
$V(\sigma_n^T,n):=\left\Vert K_n^{T}\right\Vert_1, \  n\in\mathbb{P}.$

For $n\in\mathbb{N}$, define $ \Delta t_{j,n}:=t_{j,n}-t_{j+1,n}, \qquad j=1,\ldots,n, $ where we set $t_{n+1,n}:=0$. 
Using Abel's transformation and the normalization condition \eqref{28} we obtain 
\begin{equation}\label{Abel1}
1=\sum_{j=1}^{n}t_{j,n} =\sum_{j=1}^{n-1}j\,\Delta t_{j,n} +n t_{n,n}. 
\end{equation}
Moreover, since $ D_j=jK_j-(j-1)K_{j-1}, $ another application of Abel's transformation yields 
\begin{equation}\label{Abel2}
K_n^T=\sum_{j=1}^{n-1}j\,\Delta t_{j,n}K_j+nt_{n,n}K_n. 
\end{equation}

It is also easy to check that 
\begin{equation}\label{Abel3} 
\sigma_n^T(f) =\sum_{j=0}^{n-1}\bigl(1-T_{j,n}\bigr) \widehat{f}(j)\psi_j, 
\end{equation} 
where $ T_{j,n}:=\sum_{r=1}^{j}t_{r,n}, \ j=1,\ldots,n $ and we set \(T_{0,n}:=0\).

Combining \eqref{Dbont} with Abel’s transformation for $K_{M_n}^T$ we have the following equality
\begin{eqnarray}\label{15}
K_{M_n}^T
&=&
\sum_{j=1}^{M_n}t_{j,M_n}D_j =
D_{M_n}
+\psi_{M_n-1}
\sum_{j=1}^{M_n-2}
j\,\Delta t_{M_n-j-1,M_n}\,
\overline{K_j} \\ \notag
&-&\psi_{M_n-1}
t_{1,M_n}(M_n-1)\overline{K_{M_n-1}}
\end{eqnarray}
and
\begin{eqnarray}\label{1}
\sum_{j=1}^{M_{\vert n\vert}-1} t_{j,n}D_j
&=&
\left(\sum_{j=1}^{M_{\vert n\vert}-1} t_{j,n}\right)D_{M_{\vert n\vert}}
+\psi_{M_{\vert n\vert}-1}\sum_{j=1}^{M_{\vert n\vert}-2}
j\,\Delta t_{M_{\vert n\vert}-j-1,n}\,\overline{K_j}\\ \notag
&-&\psi_{M_{\vert n\vert}-1}t_{1,n}(M_{\vert n\vert}-1)\overline{K_{M_{\vert n\vert}-1}}.
\end{eqnarray}

We also define the following maximal operators:
\begin{eqnarray*}
\sigma ^{\ast }(f)&:=&\sup_{n\in\mathbb{P}}\left\vert \sigma _{n}(f)\right\vert,\qquad
t^{\ast }(f):=\sup_{n\in\mathbb{P}}\left\vert t_{n}(f)\right\vert,\\
T^{\ast }(f)&:=&\sup_{n\in\mathbb{P}}\left\vert T_{n}(f)\right\vert,\qquad
\sigma ^{\ast,T}(f):=\sup_{n\in\mathbb{P}}\left\vert \sigma^T _{n}(f)\right\vert.
\end{eqnarray*}

A point $x$ is called an $m$-adic Lebesgue point of an integrable function $f$ if
\begin{eqnarray*}
	\lim_{n\to \infty}M_n\int_{I_n(x)}\left\vert f(t)-f(x)\right\vert d\mu(t)=0.
\end{eqnarray*}

Now we define the operator  (for details see \cite{wconv-sing} and \cite{GG})
\begin{eqnarray*}
	W_Af(x):= \sum_{s=0}^{A}M_s\sum_{r=1}^{m_s-1}\int_{I_A(x-re_s)}\left\vert f(t)-f(x)\right\vert d\mu(t).
\end{eqnarray*}
A point $x\in G_m $ is called a Vilenkin--Lebesgue point of the function $f\in L_1(G_m),$ if
$$\lim_{A \rightarrow \infty} W_Af(x)=0.$$

Since
$
I_n(x-r e_n)=I_n(x),$ $ r=1,\ldots,m_n-1,$
we obtain that
\[
M_n
\int_{I_n(x)}
|f(t)-f(x)|\,d\mu(t)
\leq
\frac{W_n f(x)}{m_n-1}
\leq
W_n f(x)\to0.
\]
Thus, every $m$-adic Vilenkin–Lebesgue point of $f\in L_1(G_m)$  is also a Lebesgue point. Since, by \cite{PTW}, $S_{M_n}f(x)\to f(x),$ as $ n\to\infty,$ at every $m$-adic-Lebesgue point, we can conclude that 
\begin{eqnarray}\label{VL-L}
S_{M_n}f(x)\to f(x),
\quad n\to\infty, \quad \text{at every Vilenkin--Lebesgue point} \quad x \quad \text{of} \quad f.
\end{eqnarray}

\section{Historical overview}

The weak type $(1,1)$ inequality of the maximal operator of Fejér means $\sigma ^{\ast }$  is due to 
Schipp \cite{Sc} for Walsh series and Pál and Simon \cite{PS}  for bounded Vilenkin
 series, from which almost everywhere convergence follows
by a standard argument (see \cite{mz}). In \cite{AT2}, it was proved that, for any set $\mathbb{E}\subset G_m$ of measure zero and every $1\le p<\infty$, there exists an explicit function $f\in L_p(G_m)$ whose Vilenkin--Fej\'er means diverge at every point of $\mathbb{E}$.

For the maximal operator of N\"orlund means, it was proved in \cite{PTW1} that if  $\{q_k:k\in\mathbb{N}\}$ is non-increasing and satisfies the condition 
\begin{equation} \label{fn0}
{1}/{Q_n}=O\left( {1}/{n}\right) ,\text{ \ \ as \ \ }
n\rightarrow\infty,
\end{equation}
then 
\begin{eqnarray}\label{weaktypenorlund}
y \mu\left\{t^{\ast}f>y \right\}\leq c\left\Vert f\right\Vert_{1}, \ \ \ f\in L_1(G_m), \ \ y>0.
\end{eqnarray}

If the sequence $\{q_k:k\in\mathbb{N}\}$ is non-decreasing, then the maximal operator associated with the corresponding N\"orlund means is of weak type $(1,1)$, that is, inequality \eqref{weaktypenorlund} holds. Consequently, the N\"orlund means converge almost everywhere for every $f\in L_1(G_m)$.

In \cite{APST}, it was proved that, for every set $\mathbb{E}\subset G_m$ of measure zero and every $1\le p<\infty$, there exists an explicit function $f\in L_p(G_m)$ whose Vilenkin--N\"orlund means, generated by monotone coefficients, diverge at every point of $\mathbb{E}$.

For the maximal operator of $T$-means, it was proved in \cite{BGN} that  if  $\{q_k:k\in\mathbb{N}\}$ is non-decreasing and satisfies the condition 
\begin{eqnarray} \label{10a} 
{q_{n-1}}/{Q_n}=O\left({1}/{n}\right) \ \ \ \text{as } \ \ \ n\rightarrow \infty,
\end{eqnarray}
then 
\begin{eqnarray}\label{weakT}
y \mu\left\{T^{\ast}(f)>y \right\}\leq c\left\Vert f\right\Vert_{1}, \ \ \ f\in L_1(G_m), \ \ y>0.
\end{eqnarray}
If ${q_k\in\mathbb{N}}$ is non-increasing, then the maximal operator of the corresponding $T$-means is of weak type $(1,1)$, i.e., \eqref{weakT} holds, and hence the $T$-means converge almost everywhere for every $f\in L_1(G_m)$. The analogous divergence result on sets of measure zero for Vilenkin--$T$ means was proved in \cite{APTY}, extending the result of \cite{APST}.

In \cite{AA} (see also \cite{PTW}, pp. 137--138), it was proved that 
if $G_m$ is a bounded Vilenkin group, $k\in\mathbb P$, and $f\in L_p(G_m)$, with $1\leq p\leq\infty$, then
\begin{eqnarray}\label{aaa}
\left\Vert \sigma_{k}(f)-f\right\Vert_{p}
\leq 5R^3\sum_{i=0}^{\vert k\vert} \frac{M_{i}}{M_{\vert k\vert}}\omega_p\left(f,\frac{1}{M_i}\right), \ \ \ \text{where} 
\ \ \ R:=\sup_{i\in\mathbb{N}}m_i.
\end{eqnarray}

Approximation by N\"orlund means and by more general related summability methods was investigated by M\'oricz and Siddiqi \cite{MS}, Areshidze and Tephnadze \cite{AT}, Blahota and K. Nagy \cite{BN1}, Hogan and Lakey \cite{Lak} and Weisz \cite{W1,NPTW,PSTW} (for Vilenkin versions, see \cite{AAB}, \cite{AAPT}, and \cite{BN-Publ}).
As special cases, M\'oricz and Siddiqi \cite{MS} obtained the earlier results given by   Blyumin \cite{Blu}, Efimov \cite{Ef},  Jastrebova \cite{J,J1}, Skvortsov \cite{Skc}, Tsutserova \cite{Tsu} and Yano \cite{Y2}. 
In particular, M\'oricz and Siddiqi \cite{MS} considered $f\in L_{p},\ 1\leq p\leq\infty$ and a nonnegative sequence $\{q_{k}: k\in\mathbb{N}\}$  such that
	\begin{equation}\label{MS1}
		\frac{n^{\gamma-1}}{Q^{\gamma}_{n}}\sum_{k=0}^{n-1}q_{k}^{\gamma}=O(1)\textrm{ for some } 1<\gamma\leq 2.
	\end{equation}
Under condition \eqref{MS1}, they proved that, if $\{q_k\}$ is non-decreasing, then
	\begin{equation}\label{MS2}
\left\|t_{n}(f)-f\right\|_{p}\leq \frac{5}{2Q_{n}}\sum_{j=0}^{|n|-1}2^{j}q_{n-2^{j}}\omega_{p}\left(f,\frac{1}{2^{j}}\right)+O\left(\omega_{p}\left(f,\frac{1}{2^{|n|}}\right)\right),
	\end{equation}
whereas, for non-increasing sequences $\{q_k:k\in\mathbb N\}$ they obtained
	\begin{equation}\label{MS4}
\|t_n f-f\|_p
\leq
\frac{5}{2Q_n}
\sum_{j=0}^{m-1}
\left(
Q_{n-2^j+1}
-
Q_{n-2^{j+1}+1}
\right)
\omega_p\left(f,\frac{1}{2^j}\right)
+
O\left(
\omega_p\left(f,\frac{1}{2^m}\right)
\right),
\end{equation}
where $n=2^m+r$, $1\leq r\leq2^m$ and $t_n$	are the N\"orlund means with respect to the Walsh system.
We note that condition \eqref{MS1} implies \eqref{Nor} and \eqref{reg}.

Areshidze and Tephnadze \cite{AT} proved inequality \eqref{MS2} with the same coefficients and without condition \eqref{MS1}.  Anakidze,  Areshidze and  Baramidze \cite{AAB} generalized this result for Vilenkin systems. In particular, they proved that if $\{q_k:k\in \mathbb{N}\}$ is non-decreasing, then,  for any $f\in L_p(G_m)$ with $1\leq p< \infty, $ 	
\begin{equation*}
\Vert t_n(f)-f\Vert_p\leq \frac{C_p}{Q_n} \sum_{i=0}^{{\vert n\vert}-1}M_iq_{n-M_i}\omega_p\left(f,\frac{1}{M_i}\right)+C_p\omega_p\left(f, \frac{1}{M_{\vert n\vert}}\right).
\end{equation*}
They also proved that if  $\{q_k:k\in \mathbb{N}\}$ is non-increasing,  then, for any $f\in L_p(G_m)$ with $1\leq p< \infty, $ 	we have the following  inequality
\begin{eqnarray}\label{12} \ \ \ \
	\Vert t_{M_n}(f)-f\Vert_p\leq C_p\sum_{s=0}^{n}\frac{M_s}{M_n}\omega _p\left(f,\frac{1}{M_s}\right)+C_p\overset{n-1}{\underset{s=0}{\sum}}\frac{(n-s)M_s}{M_{n}}\frac{q_{M_{s}}}{q_{M_{n}}}\omega _p\left(f,\frac{1}{M_s}\right).
\end{eqnarray}
Moreover, if condition \eqref{fn0} holds, then, for non-increasing sequences $\{q_k:k\in \mathbb{N}\}$,
\begin{eqnarray}\label{14} \ \ \
	\Vert t_n(f)-f\Vert_p \leq C_p\sum_{j=0}^{\vert n\vert}\frac{M_j}{M_{\vert n\vert}}\omega_p\left(f,\frac{1}{M_j}\right), \ \ \ (f\in L_p(G_m),\ \ \ 1\leq p< \infty).
\end{eqnarray}

It is known that if  $f\in L_p(G_m)$, then, for any $n\in \mathbb{N},$  
	\begin{eqnarray*} 
	\Vert T_n(f)-f\Vert_p\leq\frac{C_p}{Q_n}
	\sum_{j=0}^{|n|-1} M_{j}q_{M_j}\omega_p\left(f,\frac{1}{M_j}\right)
	+C_p\omega_p\left(f,\frac{1}{M_{|n|}}\right),	
	\end{eqnarray*}
	where $1\leq p< \infty$  and the $T$  means are generated by a non-increasing sequence $\{q_k:k\in\mathbb N\}$ (see \cite{AAPT}).
On the other hand, if the regular $T$  means are generated by a non-decreasing sequence $\{q_k:k\in\mathbb N\}$ and satisfy condition \eqref{10a},  then the inequality	
	\begin{eqnarray}\label{16}
	\Vert T_n(f)-f\Vert_p \leq C_p\sum_{j=0}^{\vert n\vert}\frac{M_j}{M_{\vert n\vert}}\omega_p\left(f,\frac{1}{M_j}\right)
	\end{eqnarray}
	holds for some constant $C_p$ depending only on $p.$
They also proved that 
\begin{eqnarray}\label{13} \ \ \ \ \ \ \ 
\Vert T_{M_n}(f)-f\Vert_p
\leq C_p\sum_{j=0}^{n}\frac{M_j}{M_n}\omega _p\left(f, \frac{1}{M_j}\right) 
+\frac{C_p}{q_{0}}\overset{n-1}{\underset{j=0}{\sum}}\frac{(n-j)q_{M_n-M_j}M_j}{M_{n}}\omega _p\left(f,\frac{1}{M_j}\right).
\end{eqnarray}

Pointwise convergence
of matrix-transform means of Vilenkin--Fourier series at
Vilenkin--Lebesgue points was studied by Goginava, Omira, Abdel-Latif
and AlKrinat \cite{GOAA}, and subsequently by Blahota and D. Nagy
\cite{BN-MN} for matrices with monotone rows.

For Walsh--Fourier series, norm and almost everywhere convergence of matrix-transform means were studied by Blahota and G'at \cite{BG}, while Blahota \cite{Bl} considered norm convergence of subsequences. More recently, Blahota and G'at \cite{BlGa} obtained norm and almost everywhere convergence and divergence results for subsequences of general matrix-transform means. In particular, they showed that the condition
\begin{equation}\label{cond-0} 
\sum_{k=1}^{n} t_{k,n}^{\,2} =O\left({1}/{n}\right) 
\end{equation} 
is sufficient for norm convergence and, along lacunary subsequences, for almost everywhere convergence of matrix-transform means of Walsh--Fourier series, without imposing monotonicity on the rows of the matrix. Blahota and Nagy \cite{BN-MN} investigated the convergence of matrix-transform means of Vilenkin--Fourier series at Vilenkin--Lebesgue points for matrices with monotone rows. In particular, for non-increasing rows they proved that the condition $t_{1,n}\to0,\ n\to\infty,$ implies $ \sigma_n^T(f;x)\to f(x) $ at every Vilenkin--Lebesgue point $x$ of $f$. For non-decreasing rows, under the additional condition 
\begin{equation}\label{cond-1} 
t_{n,n}=O\left({1}/{n}\right), 
\end{equation} 
they established the convergence $ \sigma_n^T(f;x)\to f(x) $ at every point $x$ at which the Fej\'er means converge to $f(x)$, and hence, in particular, at every Vilenkin--Lebesgue point. More recently, Szilveszter, Toledo and Weisz \cite{STW} extended pointwise convergence results for matrix-transform means to bounded multi-dimensional Vilenkin--Fourier series and established convergence at $\omega$-Vilenkin--Lebesgue points under suitable conditions on the matrix coefficients.

Blahota and K. Nagy \cite{BN-Publ} studied approximation by
matrix-transform means of Vilenkin--Fourier series. In particular, if
$f\in L_p(G_m)$, $1\leq p\leq\infty$, and
$\{t_{k,n}:1\leq k\leq n\}$ is non-increasing for every fixed $n$, then
they obtained an estimate of the form
\begin{equation}\label{11}
\left\|\sigma_n^T(f)-f\right\|_p
\leq
c\sum_{j=0}^{|n|-1}
M_j
\sum_{l=1}^{m_j-1}
t_{lM_j,n}
\omega_p\left(f,\frac{1}{M_j}\right)
+
c\sum_{k=M_{|n|}}^{n}
t_{k,n}
\omega_p\left(f,\frac{1}{M_{|n|}}\right).
\end{equation}
In particular, since
$
\sum_{k=M_{|n|}}^{n}t_{k,n}\leq1,
$
the second term in \eqref{11} is bounded by
$
c\,\omega_p\left(f,\frac{1}{M_{|n|}}\right).
$
Moreover, if the rows are non-decreasing and \eqref{cond-1}
is satisfied, 
then
\begin{equation}\label{17}
\begin{aligned}
\left\|\sigma_n^T(f)-f\right\|_p
\leq{}&
c\sum_{j=0}^{|n|-1}
M_j\left(
\sum_{l=1}^{m_j-1}
t_{(l+1)M_j-1,n}\right)
\omega_p\left(f,\frac{1}{M_j}\right)
+
O\left(
\omega_p\left(f,\frac{1}{M_{|n|}}\right)
\right).
\end{aligned}
\end{equation}

More recently, Agafonova and Volosivets \cite{AV} obtained approximation estimates for general matrix-transform means under suitable variation conditions on the matrix coefficients and an additional condition corresponding to \eqref{cond-1}. In contrast, our estimates for non-decreasing rows do not require this additional condition. In the Walsh case, related approximation results were obtained by Goginava and K. Nagy \cite{GN10} and, more recently, by Blahota and D. Nagy \cite{BN}.

We prove new convergence, divergence, and approximation results for matrix-transform means of Vilenkin--Fourier series. For bounded groups, non-decreasing rows converge along ${M_n}$ without condition \eqref{cond-1} (Theorem~\ref{BD4_5}). For arbitrary bounded or unbounded groups, every regular matrix-transform method admits divergence on any prescribed null set along any increasing sequence (Theorem~\ref{th4}). We also obtain $L_p$-approximation estimates, including a full-sequence bound involving $V(\sigma_n^T,n)$.

\section{Auxiliary results}	

We also need the following result from the paper \cite{PSTW} (see also the recent book \cite{PTW}, Theorem 2.6.4, p. 111):
\begin{lemma}\label{lemma3}
	
	Let $G_m$ be an arbitrary Vilenkin group (bounded or unbounded). Suppose that $\mathbb{E}$ is a subset of $G_m$ such that $\mu (\mathbb{E})= 0.$ Cover $\mathbb{E}$ with intervals	$(\mathbb{I}_i, i\in \mathbb{N})$ satisfying the following conditions. Each point $x \in \mathbb{E}$ belongs to infinitely many of the intervals $\mathbb{I}_k$ and
	$\sum_{i=0}^{\infty}\mu (\mathbb{I}_i)<1.$
	Set $n_0 := 0$ and choose an increasing sequence $\{n_i\}$ of integers so that
	$$\sum_{k=n_i}^{\infty}\mu (\mathbb{I}_k)<M_i^{-1} \ \  (i \in \mathbb{N}).$$
	Then, for the sets
	$\mathbb{A}_i := \bigcup_{k=n_i}^{n_{i+1}-1} \mathbb{I}_k \  (i \in \mathbb{N}) $
	there exists a sequence of  increasing integers $\{\alpha_i\}$ and a sequence 
	of Vilenkin polynomials $\{P_i, i\in\mathbb{N}\}$ such that
	$sp(P_i):=\{k\in\mathbb{N}: \widehat{P}_i(k)\neq 0 \}\subset \left[M_{\alpha_i}, M_{\alpha_{i+1}} \right),$
	that is
	$
	P_i=\sum_{k=M_{\alpha_i}}^{M_{\alpha_{i+1}}-1}c_k^i\psi_k,
	$
	\begin{equation}\label{9aa}
	\|P_i \|_p^p=\mu(\mathbb{A}_i)\leq M_i^{-1}, \ \text{for every} \ \ 1\leq p< \infty,
	\end{equation}
	and, for any $i \in \mathbb{N}$,
	\begin{equation}\label{25} 
	\left|P_i(x)\right|=1, \ \ \  x\in \mathbb{A}_i.
	\end{equation}	
\end{lemma}

We will also need the following kernel estimates for matrix-transform means:
\begin{lemma}\label{weak00} 
Let $n\in\mathbb P$, and let $G_m$ be a bounded Vilenkin group. Suppose that, for every $k\in\mathbb P$, the coefficient sequence
${t_{i,k}:1\le i\le k}$
is non-decreasing and  \eqref{28} is satisfied. Then, for every $n\in\mathbb P,$
	\begin{equation}\label{KT1}
	\left\Vert K_{M_n}^T\right\Vert_1\leq 4R^3+1,
	\end{equation}
and
\begin{equation}\label{KT2}
\left\Vert\sum_{j=1}^{M_{\left\vert n\right\vert}-1}t_{j,n}D_j\right\Vert_1\leq 4R^3+1.
\end{equation}
	
If, in addition, the coefficient sequences satisfy either condition \eqref{cond-1} or
\begin{equation}\label{MS3}
n^{\gamma-1}
\sum_{k=1}^{n}t_{k,n}^{\gamma}
=O(1),
\qquad 1<\gamma\leq2,
\end{equation}
then there exists a constant $C$ such that
\begin{equation}\label{KT3}
\left\Vert K_{n}^T\right\Vert_1\leq C<\infty.
\end{equation}
\end{lemma}

\begin{proof}
	By combining \eqref{dn40}, \eqref{fn40} and \eqref{15}
	we can conclude that
	\[
	\begin{aligned}
	\|K_{M_n}^T\|_1
	&\le
	1+4R^3\left(
	\sum_{j=1}^{M_n-2}
	j\left|\Delta t_{M_n-j-1,M_n}\right|
	+(M_n-1)t_{1,M_n}
	\right).
	\end{aligned}
	\]
	
	Since the row is non-decreasing, $
	\left|\Delta t_{M_n-j-1,M_n}\right|
	=
	t_{M_n-j,M_n}-t_{M_n-j-1,M_n},
	$
	we have that
	\[
	\begin{aligned}
	&\sum_{j=1}^{M_n-2}
	j\left|\Delta t_{M_n-j-1,M_n}\right|
	+(M_n-1)t_{1,M_n}=
	\sum_{j=1}^{M_n-2}
	j\bigl(t_{M_n-j,M_n}-t_{M_n-j-1,M_n}\bigr)
	+(M_n-1)t_{1,M_n}.
	\end{aligned}
	\]
	
	Changing the index by $r=M_n-j-1$, we obtain
	\[
	\begin{aligned}
	\sum_{j=1}^{M_n-2}
	j\bigl(t_{M_n-j,M_n}-t_{M_n-j-1,M_n}\bigr)
	&=
	\sum_{r=1}^{M_n-2}
	(M_n-r-1)\bigl(t_{r+1,M_n}-t_{r,M_n}\bigr)\\
	&=
	\sum_{r=2}^{M_n-1}t_{r,M_n}
	-(M_n-2)t_{1,M_n}.
	\end{aligned}
	\]
	Hence,
	\begin{eqnarray*}
	\begin{aligned}
	&\sum_{j=1}^{M_n-2}
	j\left|\Delta t_{M_n-j-1,M_n}\right|
	+(M_n-1)t_{1,M_n}=
	\sum_{r=1}^{M_n-1}t_{r,M_n}
	=
	1-t_{M_n,M_n}
	\le 1.
	\end{aligned}
	\end{eqnarray*}
	
	Consequently,
	\[
	\|K_{M_n}^T\|_1
	\le
	1+4R^3\bigl(1-t_{M_n,M_n}\bigr)
	\le
	4R^3+1,
	\]
	which completes the proof of \eqref{KT1}.
	
Analogously, combining \eqref{dn40}, \eqref{fn40} and \eqref{1} we can conclude that
	\[
	\begin{aligned}
	\left\|\sum_{j=1}^{M_{\left\vert n\right\vert}-1} t_{j,n}D_j\right\|_1
	&\le
	T_{M_{\left\vert n\right\vert}-1,n}
	+4R^3\left(
	\sum_{j=1}^{M_{\left\vert n\right\vert}-2}
	j\bigl(t_{M_{\left\vert n\right\vert}-j,n}
	-t_{M_{\left\vert n\right\vert}-j-1,n}\bigr)
	+(M_{\left\vert n\right\vert}-1)t_{1,n}
	\right)
	\\
	&=
	T_{M_{\left\vert n\right\vert}-1,n}
	+4R^3T_{M_{\left\vert n\right\vert}-1,n}
	=
	(4R^3+1)T_{M_{\left\vert n\right\vert}-1,n}
	\le 4R^3+1.
	\end{aligned}
	\]
	This proves estimate \eqref{KT2}.
	
	Since the row is non-decreasing, by using \eqref{Abel1} we obtain that
	$
	|\Delta t_{j,n}|
	=
	t_{j+1,n}-t_{j,n}
	$
	and 
	\[
	\sum_{j=1}^{n-1}
	j|\Delta t_{j,n}|
	=
	\sum_{j=1}^{n-1}
	j(t_{j+1,n}-t_{j,n})
	=
	nt_{n,n}-1.
	\]
	By combining \eqref{fn40} and \eqref{Abel2} we find that
	\[\begin{aligned}
	\|K_n^T\|_1\leq 4R^3\left(\sum_{j=1}^{n-1}j|\Delta t_{j,n}|+nt_{n,n}\right)=4R^3(2nt_{n,n}-1)\leq C.
	\end{aligned}\]
	Thus, estimate \eqref{KT3} holds under condition \eqref{cond-1}.
	
	By the Sidon-type inequality for bounded Vilenkin systems
	(see Daly and Fridli \cite[Lemma~8]{DF}), for $1<\gamma\leq2$ there exists a constant
	$C_{\gamma,R}>0$ such that
	\[
	\left\|
	\sum_{k=1}^{n}a_kD_k
	\right\|_1
	\leq
	C_{\gamma,R}
	n^{1-1/\gamma}
	\left(
	\sum_{k=1}^{n}|a_k|^\gamma
	\right)^{1/\gamma}.
	\]
	
	Applying the preceding inequality with $a_k=t_{k,n}$  
	we can conclude that
	\[
	\begin{aligned}
	\|K_n^T\|_1\leq
	C_{\gamma,R}
	n^{1-1/\gamma}
	\left(
	\sum_{k=1}^{n}t_{k,n}^{\gamma}
	\right)^{1/\gamma}=
	C_{\gamma,R}
	\left(
	n^{\gamma-1}
	\sum_{k=1}^{n}t_{k,n}^{\gamma}
	\right)^{1/\gamma}.
	\end{aligned}
	\]
	Hence, condition \eqref{MS3} also implies \eqref{KT3}. This completes the proof of the lemma.
\end{proof}

\section{Pointwise and almost everywhere convergence of matrix-transform means}

The next theorem removes condition \eqref{cond-1} along the subsequence $\{M_n\}$:

\begin{theorem}\label{BD4_5}
	Let $f\in L_1(G_m)$, where $G_m$ is a bounded Vilenkin group. Suppose that, for every fixed
	$n\in\mathbb P$, the row $\{t_{i,n}:1\leq i\leq n\}$ is
	non-decreasing and \eqref{28} is satisfied. Then
	\[
	\sigma^T_{M_n}(f;x)\to f(x)
	\]
	at every Vilenkin--Lebesgue point $x$ of $f$ as $n\to\infty$.
\end{theorem}
\begin{proof}
Under the assumptions of the theorem, the matrix-transform method is
automatically regular.
Indeed, for every fixed $i\in\mathbb{P}$ and $n\geq i$, since the sequence
$\{t_{j,n}:1\leq j\leq n\}$ is non-decreasing and \eqref{28} holds, we have
$
(n-i+1)t_{i,n}
\leq \sum_{j=i}^{n}t_{j,n}
\leq 1.
$
Hence,
$t_{i,n}\leq \frac{1}{n-i+1}\to 0,$
 as $n\to\infty.$
Therefore, condition \eqref{29} is satisfied. Moreover,
\begin{eqnarray}\label{22aa}
t_{1,M_{n}}(M_{n}-1)\leq \sum_{i=1}^{M_n-1}t_{i,M_{n}}\leq \sum_{i=1}^{M_n}t_{i,M_{n}}=1,
\end{eqnarray}
and 
\begin{eqnarray}\label{22a}
&&\sum_{j=1}^{M_n-2}
j\left|\Delta t_{M_n-j-1,M_n}\right|\\ \notag
&=&\sum_{j=1}^{M_n-2}j (t_{M_n-j,M_n}-t_{M_n-j-1,M_n})=
\sum_{j=1}^{M_n-2}
(M_n-j-1)
\left(t_{j+1,M_n}-t_{j,M_n}\right)\\ \notag
&=&
\sum_{j=2}^{M_n-1} t_{j,M_n}
-(M_n-2)t_{1,M_n}=
1-t_{M_n,M_n}-(M_n-1)t_{1,M_n}\leq 1.
\end{eqnarray}

Now fix a Vilenkin--Lebesgue point $x$ of $f$. Define
$
g_j(u):=[f(x-u)-f(x)]\overline{K_j(u)}.
$
Then
\begin{align}\label{12a}
&\sigma_{M_n}^T(f;x)-f(x)
= S_{M_n}f(x)-f(x) \notag\\ &+\sum_{j=1}^{M_n-2} j\,\Delta t_{M_n-j-1,M_n} \int_{G_m} \bigl(f(x-u)-f(x)\bigr) \psi_{M_n-1}(u)\overline{K_j(u)} \,d\mu(u) \notag\\ &-t_{1,M_n}(M_n-1) \int_{G_m} \bigl(f(x-u)-f(x)\bigr) \psi_{M_n-1}(u)\overline{K_{M_n-1}(u)} \,d\mu(u)\notag\\
={}& S_{M_n}f(x)-f(x) \notag+\sum_{j=1}^{M_n-2}
j\,\Delta t_{M_n-j-1,M_n}
\int_{G_m}
g_j(u)\psi_{M_n-1}(u)
\,d\mu(u) \notag\\
&-t_{1,M_n}(M_n-1)
\int_{G_m}
g_{M_n-1}(u)\psi_{M_n-1}(u)
\,d\mu(u)
\notag\\ 
={}& \left(S_{M_n}f(x)-f(x)\right)+\sum_{j=1}^{J-1}j\,\Delta t_{M_n-j-1,M_n}\int_{G_m}g_j(u)\psi_{M_n-1}(u)
\,d\mu(u) \notag\\
&+\sum_{j=J}^{M_n-2}j\,\Delta t_{M_n-j-1,M_n}\int_{G_m}g_j(u)\psi_{M_n-1}(u)\,d\mu(u) \notag\\
&-t_{1,M_n}(M_n-1)
\int_{G_m}
g_{M_n-1}(u)\psi_{M_n-1}(u)
\,d\mu(u)\notag\\
&=I_1+I_2+I_3-I_4.
\end{align}

Since $x$ is a Vilenkin--Lebesgue point, \ref{VL-L} gives
$
S_{M_n}f(x)\to f(x),
\quad n\to\infty.
$
Hence
\begin{align}\label{AE1}
|I_1|\to0.
\end{align}

For each fixed $j$, $g_j\in L_1(G_m)$. Hence the Riemann--Lebesgue theorem gives
\[
R_{n,j}(x)
:=
\int_{G_m}
g_j(u)\psi_{M_n-1}(u)\,d\mu(u)
\rightarrow 0
\qquad (n\to\infty).
\]

Since only finitely many $j<J$ occur,
$
\max_{1\le j<J}|R_{n,j}(x)|
\rightarrow 0.
$
Hence, for sufficiently large $n,$ we get 
$
\max_{1\le j<J}|R_{n,j}(x)|
\leq \varepsilon.
$
By applying \eqref{22a} we obtain that
\begin{align}\label{AE2}
\vert I_2\vert&\leq\sum_{j=1}^{J-1}j\left\vert\Delta t_{M_n-j-1,M_n}\right\vert\left\vert\int_{G_m}g_j(u)\psi_{M_n-1}(u)
\,d\mu(u)\right\vert \\ \notag
&\leq\sum_{j=1}^{J-1}j\left\vert\Delta t_{M_n-j-1,M_n}\right\vert\varepsilon \leq\varepsilon\sum_{j=1}^{M_n-2}j\left\vert\Delta t_{M_n-j-1,M_n}\right\vert \leq\varepsilon.
\end{align}

Let
\[
E_j(x):=
\int_{G_m}
\left|f(x-u)-f(x)\right|
\left|\overline{K_j(u)}\right|d\mu(u)=\int_{G_m}
\left|f(x-u)-f(x)\right|
\left|K_j(u)\right|
\,d\mu(u).
\]

By using estimates in \eqref{fn41}	we obtain that
\[
\begin{aligned}
E_n(x)
&\leq
\frac{c}{n}
\sum_{A=0}^{|n|}
M_A
\int_{G_m}
|f(x-u)-f(x)|\,|K_{M_A}(u)|\,d\mu(u)\leq
\frac{c}{n}
\sum_{A=0}^{|n|}
M_A W_Af(x).
\end{aligned}
\]
Since $x$ is a Vilenkin--Lebesgue point,
$
W_Af(x)\to0,
\ A\to\infty.
$
Hence, by splitting the last sum into a fixed initial part and a tail,
and using
$
\sum_{A=0}^{|n|}M_A\leq cM_{|n|}\leq cn,
$
we conclude that
$
E_n(x)\to0,
\ n\to\infty.
$
For every $\varepsilon>0$, choose $J$ such that
$
E_j(x)<\varepsilon,
\ j\ge J.
$
Then
\begin{eqnarray*}
\sup_{J\le j<M_n}|E_{j}(x)|<\varepsilon.
\end{eqnarray*}

By using \eqref{22aa} and  \eqref{22a} we get that
\begin{align}\label{AE3}
\left\vert I_3\right\vert+\left\vert I_4 \right\vert&\leq
\sum_{j=J}^{M_n-2}
j\left|\Delta t_{M_n-j-1,M_n}\right|E_{j}(x)
+
t_{1,M_n}(M_n-1)E_{M_n-1}(x)\\ \notag
&\leq\varepsilon\left(\sum_{j=J}^{M_n-2}
j\left|\Delta t_{M_n-j-1,M_n}\right|
+
t_{1,M_n}(M_n-1)\right)\\ \notag
&\leq\varepsilon\left(\sum_{j=1}^{M_n-2}
j\left|\Delta t_{M_n-j-1,M_n}\right|
+
t_{1,M_n}(M_n-1)\right)\leq 2\varepsilon.
\end{align}

Since $\varepsilon>0$ is arbitrary, by combining \eqref{12a}, \eqref{AE1}, \eqref{AE2} and \eqref{AE3} we obtain that, at every Vilenkin--Lebesgue point,
$
\sigma^T_{M_n}(f;x)\to f(x),
\ n\to\infty.
$
The proof is complete.

\end{proof}

In the next theorem, we establish a complementary divergence result for arbitrary regular matrix-transform methods along any prescribed increasing sequence:
\begin{theorem}\label{th4}
	Let $G_m$ be an arbitrary Vilenkin group (bounded or unbounded), let $\{\sigma_n^T\}$ be a regular matrix-transform summability method satisfying \eqref{28}, and let $\{n_k\}_{k\in\mathbb N}$ be an arbitrary increasing sequence of positive integers. If $\mathbb E\subset G_m$ has measure zero, then there exists
	$
	f\in\bigcap_{1\leq p<\infty}L_p(G_m)
	$
	such that the sequence
	$
	\{\sigma_{n_k}^T(f;x)\}_{k\in\mathbb N}
	$
	diverges at every point $x\in\mathbb E$.
\end{theorem}
\begin{proof}
	Let $\{\mathbb A_i\}_{i\in\mathbb N}$,
	$\{\alpha_i\}_{i\in\mathbb N}$ and
	$\{P_i\}_{i\in\mathbb N}$ be as in Lemma~\ref{lemma3}. Thus,
	$
	\operatorname{sp}(P_i)
	\subset
	[M_{\alpha_i},M_{\alpha_{i+1}}),
	$
	and
	$
	\|P_i\|_p^p
	=
	\mu(\mathbb A_i)
	\leq M_i^{-1},$ $ 1\leq p<\infty,
	$
	while
	$
	|P_i(x)|=1,$ $ x\in\mathbb A_i.
	$
	Moreover, every $x\in\mathbb E$ belongs to infinitely many of the
	sets $\mathbb A_i$.
	
	Recall that
	$
	T_{r,n}:=\sum_{s=1}^{r}t_{s,n}.
	$
	Since the matrix-transform summability method is regular, by
	\eqref{29}, we get that
	$t_{s,n}\to0$ as $ n\to\infty$
	for every fixed $s\in\mathbb P$. Hence, for every fixed
	$L\in\mathbb P$,
	\begin{equation}\label{TL0}
	T_{L,n}
	=
	\sum_{s=1}^{L}t_{s,n}
	\to0,
	\qquad n\to\infty.
	\end{equation}
	
	We now construct recursively a subsequence
	$\{n_{k_i}\}_{i\in\mathbb N}$ of the prescribed sequence
	$\{n_k\}_{k\in\mathbb N}$ and an increasing sequence of integers
	$\{\beta_i\}_{i\in\mathbb N}$.
	
	Choose $k_0\in\mathbb N$. After $k_i$ has been chosen, choose
	$\beta_i$ sufficiently large so that
	$
	\beta_i>\alpha_{i+1},$
	$
	M_{\beta_i}>n_{k_i},
	$
	and, for $i\geq1$, also
	$
	\beta_i>\beta_{i-1}.
	$
	Set
	$
	\xi_i:=
	M_{\beta_i}+M_{\alpha_{i+1}}-1.
	$
	Since $\{n_k\}$ is increasing and tends to infinity, and since
	\eqref{TL0} holds for the fixed integer $\xi_i$, we may choose
	$k_{i+1}>k_i$ sufficiently large so that
	\begin{equation}\label{nkchoice}
	n_{k_{i+1}}>\xi_i
	\end{equation}
	and
	\begin{equation}\label{Tsmall}
	T_{\xi_i,n_{k_{i+1}}}
	<{1}/{\xi_i^3}.
	\end{equation}
	
	Let
	$
	g_i:=\psi_{M_{\beta_i}}P_i
	$
	and define
	\begin{equation}\label{construction2}
	f:=
	\sum_{i=0}^{\infty}g_i.
	\end{equation}
	Since
	$\|g_i\|_p=\|P_i\|_p\leq M_i^{-1/p}\leq 2^{-i/p},$
	we obtain
	$
	\sum_{i=0}^{\infty}
	\|g_i\|_p
	<\infty
	$
	for every $1\leq p<\infty$. Hence the series in
	\eqref{construction2} converges in $L_p(G_m)$ for every
	$1\leq p<\infty$, and 
	$
	f\in
	\bigcap_{1\leq p<\infty}L_p(G_m).
	$
Since $\beta_i>\alpha_{i+1}$, for
$M_{\alpha_i}\leq r<M_{\alpha_{i+1}}$ we have
$r<M_{\beta_i}$, and hence
$
\psi_{M_{\beta_i}}\psi_r
=
\psi_{M_{\beta_i}+r}.
$
Therefore,
$
\operatorname{sp}(g_i)
\subset
\left[
M_{\beta_i}+M_{\alpha_i},
M_{\beta_i}+M_{\alpha_{i+1}}
\right).
$
Moreover, since $\{\beta_i\}$ is strictly increasing,
$
M_{\beta_i}+M_{\alpha_{i+1}}
<
M_{\beta_{i+1}}+M_{\alpha_{i+1}},
$
and hence these spectral blocks are pairwise disjoint.
Since the series in \eqref{construction2} converges in $L_1(G_m)$, we have
\[
\widehat f(k)
=
\begin{cases}
\widehat g_l(k),
&
M_{\beta_l}+M_{\alpha_l}
\leq k
<
M_{\beta_l}+M_{\alpha_{l+1}},
\quad l\in\mathbb N,
\\[2mm]
0,
&
\text{otherwise}.
\end{cases}
\]
Moreover,
$
\widehat g_l(k)
=
c_{k-M_{\beta_l}}^l
$
on the corresponding spectral block. Hence
\[
\widehat f(k)
=
\begin{cases}
c_{k-M_{\beta_l}}^l,
&
M_{\beta_l}+M_{\alpha_l}
\leq k
<
M_{\beta_l}+M_{\alpha_{l+1}},
\quad l\in\mathbb N,
\\[2mm]
0,
&
\text{otherwise}.
\end{cases}
\]

	By the choice of $\beta_i$ and \eqref{nkchoice}, we get that
	$
	n_{k_i}
	<
	M_{\beta_i}+M_{\alpha_i}$
	and
	$
	M_{\beta_i}+M_{\alpha_{i+1}}-1
	=
	\xi_i
	<
	n_{k_{i+1}}.
	$
	Hence,
	\begin{equation}\label{blocklocation}
	\operatorname{sp}
	\left(
	g_i
	\right)
	\subset
	(n_{k_i},n_{k_{i+1}}).
	\end{equation}
	By combining \eqref{Abel3} and \eqref{blocklocation} we can conclude that
	\begin{align}\label{sigma}
	&\sigma_{n_{k_{i+1}}}^T(f;x)
	-
	\sigma_{n_{k_i}}^T(f;x)\\ \notag
	&=
	g_i(x)
	-
	\sum_{r\in\operatorname{sp}(g_i)}
	T_{r,n_{k_{i+1}}}
	\widehat g_i(r)\psi_r(x)+
	\sum_{l=0}^{i-1}
	\sum_{r\in\operatorname{sp}(g_l)}
	\left(
	T_{r,n_{k_i}}
	-
	T_{r,n_{k_{i+1}}}
	\right)
	\widehat g_l(r)\psi_r(x)
	\\ \notag
	&=:
	\psi_{M_{\beta_i}}P_i(x)
	-
	\sum_{k=M_{\beta_i}+M_{\alpha_i}}^{M_{\beta_i}+M_{\alpha_{i+1}}-1}
	T_{k,n_{k_{i+1}}}
	c^i_{k-M_{\beta_i}}\psi_k(x)
	\\ \notag
	&+
	\sum_{l=0}^{i-1}
	\sum_{k=M_{\beta_l}+M_{\alpha_l}}^{M_{\beta_l}+M_{\alpha_{l+1}}-1}
	\left(
	T_{k,n_{k_i}}
	-
	T_{k,n_{k_{i+1}}}
	\right)
	c^l_{k-M_{\beta_l}}\psi_k(x)
	\\ \notag
	&:= 
	II_1-II_2+II_3.
\end{align}
	
	We estimate $II_2$ and $II_3$. For
	$M_{\alpha_i}\leq r<M_{\alpha_{i+1}}$, by
	\eqref{9aa},
	\[
	|c_r^i|
	=
	|\widehat P_i(r)|
	\leq
	\|P_i\|_1
	=
	\mu(\mathbb A_i)
	\leq
	M_i^{-1}
	\leq1.
	\]
	Consequently, for
	$
	M_{\beta_i}+M_{\alpha_i}
	\leq k
	<
	M_{\beta_i}+M_{\alpha_{i+1}},
	$
	we have
	$
	|\widehat g_i(k)|
	=
	\left|c^i_{k-M_{\beta_i}}\right|
	\leq1.
	$
	Therefore,
	$
	|\widehat g_i(r)|\leq1.
	$
	Hence, using \eqref{Tsmall} with estimate
	$
	k\leq
	M_{\beta_i}+M_{\alpha_{i+1}}-1
	=
	\xi_i,
	$
	we obtain that
	\begin{eqnarray}\label{IIa}
	|II_2|
	&\leq&
	\sum_{k=M_{\beta_i}+M_{\alpha_i}}^
	{M_{\beta_i}+M_{\alpha_{i+1}}-1}
	T_{k,n_{k_{i+1}}}
	\left|c^i_{k-M_{\beta_i}}\right|
	\\ \notag
	&\leq&
	\sum_{k=0}^{\xi_i}
	T_{k,n_{k_{i+1}}}\leq
	(\xi_i+1)T_{\xi_i,n_{k_{i+1}}}
	\leq
	\frac{\xi_i+1}{\xi_i^3}
	\leq
	\frac{2}{\xi_i^2}
	\rightarrow0,
	\qquad i\to\infty.
	\end{eqnarray}

	For $II_3$, since all the spectra of
	$g_0,\ldots,g_{i-1}$ are contained in
	$\{0,\ldots,\xi_{i-1}\}$, we have
	\begin{eqnarray*}
	|II_3|
	&\leq&
	\sum_{l=0}^{i-1}
	\sum_{k=M_{\beta_l}+M_{\alpha_l}}^
	{M_{\beta_l}+M_{\alpha_{l+1}}-1}
	\left|
	T_{k,n_{k_i}}
	-
	T_{k,n_{k_{i+1}}}
	\right|
	\left|c^l_{k-M_{\beta_l}}\right|
	\\ \notag
	&\leq&
	\sum_{k=0}^{\xi_{i-1}}
	\left(
	T_{k,n_{k_i}}
	+
	T_{k,n_{k_{i+1}}}
	\right)
	\leq
	(\xi_{i-1}+1)
	\left(
	T_{\xi_{i-1},n_{k_i}}
	+
	T_{\xi_{i-1},n_{k_{i+1}}}
	\right).
	\end{eqnarray*}
	By using \eqref{Tsmall} we get that
	$
	T_{\xi_{i-1},n_{k_i}}
	<
	{1}/{\xi_{i-1}^3}.
	$
	Moreover, since $\xi_{i-1}\leq \xi_i$,
	\[
	T_{\xi_{i-1},n_{k_{i+1}}}
	\leq
	T_{\xi_i,n_{k_{i+1}}}
	<
	\frac{1}{\xi_i^3}
	\leq
	\frac{1}{\xi_{i-1}^3}.
	\]
	Consequently,
    \begin{eqnarray}\label{AA}
	|II_3|
	\leq
	{2}/{\xi_{i-1}^2}
	\rightarrow0.
	\end{eqnarray}
	
	Finally, every $x\in\mathbb E$ belongs to infinitely many of the
	sets $\mathbb A_i$. For each such $i$, by \eqref{25},
	\begin{eqnarray}\label{III}
	|II_1|=|g_i(x)|=|\psi_{M_{\beta_i}}(x)|\,|P_i(x)|=1.
	\end{eqnarray}
	
	By combining \eqref{sigma}, \eqref{IIa}, \eqref{AA} and \eqref{III}, for every $x\in\mathbb E$ and for all sufficiently
	large $i$ belonging to this infinite set,
	\[
	\begin{aligned}
	\left|
	\sigma_{n_{k_{i+1}}}^T(f;x)
	-
	\sigma_{n_{k_i}}^T(f;x)
	\right|\geq
	|II_1|-|II_2|-|II_3|
	\geq
	1-\frac{2}{\xi_i^2}
	-\frac{2}{\xi_{i-1}^2}
	\geq
	\frac12.
	\end{aligned}
	\]
	Thus the sequence
	$
	\{\sigma_{n_{k_i}}^T(f;x)\}_{i\in\mathbb N}
	$
	is not a Cauchy sequence for any $x\in\mathbb E$.
	Since $\{n_{k_i}\}$ is a subsequence of the arbitrarily prescribed
	sequence $\{n_k\}$, it follows that
	$
	\{\sigma_{n_k}^T(f;x)\}_{k\in\mathbb N}
	$
	diverges at every point $x\in\mathbb E$.
	The proof is complete.
\end{proof}


\section{On the rate of approximation}

In the next theorem, in contrast to the approximation result of Agafonova and Volosivets \cite{AV}, no condition corresponding to \eqref{cond-1} is required. 
\begin{theorem}\label{BD4_2} 
	Let $G_m$ be a bounded Vilenkin group, and let $f\in L_p(G_m)$ with $1\le p\le\infty$. If  $\{t_{i,k}: 1\leq i\leq k\}$ is non-decreasing for every fixed $k\in\mathbb{P}$ and \eqref{28} is satisfied,
then 		
	\begin{eqnarray} \label{3.1} 
	&&\Vert \sigma^{T}_{M_{n}}(f)-f\Vert_p\leq (4R^3+1) \sum_{l=0}^{n-1}\left(T_{M_{l+1}-1,M_n}-T_{M_{l}-1,M_n} \right)\omega_p\left(f,1/M_l\right)+\omega_p\left(f,1/M_n\right).
	\end{eqnarray}
\end{theorem}
\begin{proof}	
	It is evident that
	\begin{eqnarray}\label{3.2}	 
	f - \sigma^{T}_{M_{n}}(f)&=&f -\sum_{j=1}^{M_n}t_{j,M_n}S_jf \\ \notag
	&=& \sum_{l=0}^{n-1}\sum_{j=M_l}^{M_{l+1}-1}t_{j,M_n}\left(f-S_jf\right)
	+t_{M_n,M_n}\left(f-S_{M_n}f\right)=IV_1+IV_2.
	\end{eqnarray}
	For $IV_2$ we have
		\begin{eqnarray}	\label{47}
	\left\Vert IV_2\right\Vert_p&\leq& \omega_p\left(f,1/M_n\right).
	\end{eqnarray}
	 
	We can rewrite $IV_1$ as
	\begin{eqnarray}\label{3.0}	
	IV_1&=&\left(\sum_{l=0}^{n-1}\sum_{j=M_l}^{M_{l+1}-1}t_{j,M_n}\left(f-S_{M_{l+1}}f\right)\right)  +\left(\sum_{l=0}^{n-1}\sum_{j=M_l}^{M_{l+1}-1}t_{j,M_n}\left(S_{M_{l+1}}f-S_jf\right)\right)
	\\ \notag
	&=&\left(\sum_{l=0}^{n-1} \left(\sum_{j=M_l}^{M_{l+1}-1}t_{j,M_n}\right)\left(f-S_{M_{l+1}}f\right)
	\right)+\sum_{l=0}^{n-1}\sum_{j=M_l}^{M_{l+1}-1}t_{j,M_n}
	\left(S_{M_{l+1}}f-S_jf\right)\\ \notag
	&=:& IV_1^1 + IV_1^2.
	\end{eqnarray}
	
	Since 
	\begin{eqnarray}\label{3.30}	
	\omega_p\left(f,1/M_{l+1}\right)\leq \omega_p\left(f,1/M_l\right), \ \ \text{for any} \  \ l\in\mathbb{N},
	\end{eqnarray}
	we obtain that 
	\begin{eqnarray}\label{3.3}		
	\left\Vert IV_1^1\right\Vert_p&\leq& \left( \sum_{l=0}^{n-1}\left(\sum_{j=M_l}^{M_{l+1}-1}t_{j,M_n}\right)\omega_p\left(f,1/M_{l+1}\right)\right)\\ \notag
	&\leq&\sum_{l=0}^{n-1}\left(T_{M_{l+1}-1,M_n}-T_{M_{l}-1,M_n} \right)\omega_p\left(f,1/M_l\right).
	\end{eqnarray}
	
	By applying \eqref{Dbont} we obtain that	
	\begin{eqnarray*}
		IV_1^2&=&\sum_{l=0}^{n-1}\sum_{j=M_l}^{M_{l+1}-1}t_{j,M_n}\left(S_{M_{l+1}}\left(f-S_{M_l}f\right)-S_j\left(f-S_{M_l}f\right)\right)\\
		&=&\sum_{l=0}^{n-1}\sum_{j=M_l}^{M_{l+1}-1}\left(t_{j,M_n}\left(D_{M_{l+1}}-D_j\right)\right)\ast\left(f-S_{M_l}f\right)\\
		&=&\sum_{l=0}^{n-1}\sum_{j=M_{l}}^{M_{l+1}-1}\left(t_{j,M_n}\psi_{M_{l+1}-1}\overline{D}_{M_{l+1}-j}\right)\ast\left(f-S_{M_l}f\right)\\
		&=&\sum_{l=0}^{n-1}\sum_{j=1}^{M_{l+1}-M_l}\left(t_{M_{l+1}-j,M_n}\psi_{M_{l+1}-1}\overline{D}_{j}\right)\ast\left(f-S_{M_l}f\right).
	\end{eqnarray*}
	Hence,		
	\begin{eqnarray*}
		\left\Vert\sum_{j=1}^{M_{l+1}-M_l}\left(t_{M_{l+1}-j,M_n}\psi_{M_{l+1}-1}\overline{D}_{j}\right)\ast
		\left(f-S_{M_l}f\right)\right\Vert_p&\leq& \left\Vert\sum_{j=1}^{M_{l+1}-M_l}t_{M_{l+1}-j,M_n}\overline{D}_{j}\right\Vert_1\left\Vert f-S_{M_l}f\right\Vert_p.
	\end{eqnarray*}
	By Abel's transformation, we obtain that
\begin{eqnarray*}
	\sum_{j=1}^{M_{l+1}-M_l}
	t_{M_{l+1}-j,M_n}D_j
	&=&
	\sum_{j=1}^{M_{l+1}-M_l-1}
	\left(
	t_{M_{l+1}-j,M_n}
	-
	t_{M_{l+1}-j-1,M_n}
	\right)jK_j
	\\
	&&+
	\left(M_{l+1}-M_l\right)
	t_{M_l,M_n}
	K_{M_{l+1}-M_l}.
\end{eqnarray*}
Applying \eqref{fn40}, we obtain that	
\begin{eqnarray*}
	&&\left\Vert
	\sum_{j=1}^{M_{l+1}-M_l}
	t_{M_{l+1}-j,M_n}\overline{D}_{j}
	\right\Vert_1
	=
	\left\Vert
	\sum_{j=1}^{M_{l+1}-M_l}
	t_{M_{l+1}-j,M_n}D_{j}
	\right\Vert_1
	\\
	&\leq&
	4R^3
	\left[
	\sum_{j=1}^{M_{l+1}-M_l-1}
	\left(
	t_{M_{l+1}-j,M_n}
	-
	t_{M_{l+1}-j-1,M_n}
	\right)j
	+
	\left(M_{l+1}-M_l\right)t_{M_l,M_n}
	\right]
	\\
	&=&
	4R^3
	\sum_{j=1}^{M_{l+1}-M_l}
	t_{M_{l+1}-j,M_n}
	=
	4R^3
	\sum_{j=M_l}^{M_{l+1}-1}
	t_{j,M_n}.
\end{eqnarray*}	
	It follows that
	\begin{eqnarray}\label{3.4}
	\left\Vert IV_1^2 \right\Vert_p&\leq& 4R^3\sum_{l=0}^{n-1}\left(\sum_{j=M_l}^{M_{l+1}-1}t_{j,M_n}\right)\omega_p\left(f,1/M_l\right)\\ \notag
	&=&4R^3\sum_{l=0}^{n-1}\left(T_{M_{l+1}-1,M_n}-T_{M_{l}-1,M_n} \right)\omega_p\left(f,1/M_l\right).
	\end{eqnarray}
	By using Minkowski's inequality in \eqref{3.2}, \eqref{47}, \eqref{3.0},	 \eqref{3.3} and \eqref{3.4} we obtain \eqref{3.1}, which completes the proof.
\end{proof}


\begin{theorem}\label{BD4_3} 
Let $G_m$ be a bounded Vilenkin group, let $f\in L_p(G_m)$ with $1\le p\le\infty$, and let $n\in\mathbb P$. Set $N:=|n|$, so that
$
M_N\leq n<M_{N+1}.
$
Suppose that, for every fixed $k\in\mathbb P$, the row
$\{t_{i,k}:1\leq i\leq k\}$ is non-decreasing and \eqref{28} is
satisfied. Then	
	\begin{eqnarray*}
\Vert \sigma^T_n f-f\Vert_p &\leq&(4R^3+1)\sum_{l=0}^{N-1}\left(T_{M_{l+1}-1,n}-T_{M_l-1,n} \right)\omega_p\left(f,1/M_l\right)\\
&+&\left(V(\sigma_{n}^T,n)+4R^3+3\right)\omega_p\left(f,1/M_{N}\right).
	\end{eqnarray*}
\end{theorem}
\begin{proof}	
It is evident that
	\begin{eqnarray*}	
f - \sigma^T_nf &=& f-\sum_{j=1}^{n}t_{j,n}S_jf=\sum_{l=0}^{N-1}\sum_{j=M_l}^{M_{l+1}-1}t_{j,n}\left(f-S_jf\right)+\sum_{j=M_N}^{n}t_{j,n}\left(f-S_jf\right)\\
&=&\sum_{l=0}^{N-1}\sum_{j=M_l}^{M_{l+1}-1}t_{j,n}\left(f-S_{M_{l+1}}f\right)
+\sum_{l=0}^{N-1}\sum_{j=M_l}^{M_{l+1}-1}t_{j,n}\left(S_{M_{l+1}}f-S_jf\right)\\
&+&\sum_{j=M_N}^{n}t_{j,n}\left(f-S_{M_N}f\right)+\sum_{j=M_N}^{n}t_{j,n}\left(S_jf-S_{M_N}f\right)\\	
		&=& V_1 + V_2+V_3+V_4 .
	\end{eqnarray*}
	
As in \eqref{3.3}, using \eqref{3.30}, we obtain that
\begin{eqnarray*}	
\left\Vert V_1\right\Vert_p&\leq&  \sum_{l=0}^{N-1} \left(\sum_{j=M_l}^{M_{l+1}-1}t_{j,n}\right)\omega_p\left(f,1/M_l\right)\leq \sum_{l=0}^{N-1}\left(T_{M_{l+1}-1,n}-T_{M_l-1,n} \right)\omega_p\left(f,1/M_l\right).
\end{eqnarray*}

If we argue as in \eqref{3.4}, by applying Abel's transformation and \eqref{fn40}, we find that
\begin{eqnarray*}
	\left\Vert
	\sum_{j=1}^{M_{l+1}-M_l}
	t_{M_{l+1}-j,n}D_j
	\right\Vert_1
	&\leq&
	4R^3
	\left[
	\sum_{j=1}^{M_{l+1}-M_l-1}
	\left(
	t_{M_{l+1}-j,n}
	-
	t_{M_{l+1}-j-1,n}
	\right)j
	+
	\left(M_{l+1}-M_l\right)t_{M_l,n}
	\right]
	\\
	&=&
	4R^3
	\sum_{j=1}^{M_{l+1}-M_l}
	t_{M_{l+1}-j,n}=
	4R^3
	\sum_{j=M_l}^{M_{l+1}-1}t_{j,n}.
\end{eqnarray*}

	It follows that
	\begin{eqnarray*}
\left\Vert V_2 \right\Vert_p&\leq& 4R^3\sum_{l=0}^{N-1}\left(\sum_{j=M_l}^{M_{l+1}-1}t_{j,n}\right)\omega_p\left(f,1/M_l\right)\\
&=&4R^3 \sum_{l=0}^{N-1}\left(T_{M_{l+1}-1,n}-T_{M_{l}-1,n} \right)\omega_p\left(f,1/M_l\right).
	\end{eqnarray*}

On the other hand, we can write
\begin{eqnarray*}	
	\left\Vert V_3\right\Vert_p&\leq& \left(\sum_{j=M_N}^{n} t_{j,n}\right)\omega_p\left(f,1/M_N\right)\leq\omega_p\left(f,1/M_N\right).
\end{eqnarray*}

For $V_4$ we get that
	\begin{eqnarray*}
		V_4&=&\sum_{j=M_N}^{n}t_{j,n}\left(S_j\left(S_{M_{N+1}}f-S_{M_{N}}f\right)-S_{M_N}\left(S_{M_{N+1}}f-S_{M_{N}}f\right)\right)\\
		&=&\sum_{j=M_N}^{n}t_{j,n}\left(D_j-D_{M_N}\right)\ast\left(S_{M_{N+1}}f-S_{M_N}f\right)
	\end{eqnarray*}
By using estimate \eqref{KT2} of Lemma~\ref{weak00}, we get that
\begin{eqnarray*}
&&\left\Vert\sum_{j=M_N}^{n}t_{j,n}\left(D_j-D_{M_N}\right)\right\Vert_1
\leq\left\Vert\sum_{j=M_N}^{n}t_{j,n}D_j\right\Vert_1+\sum_{j=M_N}^{n}t_{j,n}\leq\left\Vert\sum_{j=M_N}^{n}t_{j,n}D_j\right\Vert_1+1\\
&=&\left\Vert K_n^T-\sum_{j=1}^{M_N-1}t_{j,n}D_j\right\Vert_1+1\leq\left\Vert K_n^T\right\Vert_1+
\left\Vert\sum_{j=1}^{M_N-1}t_{j,n}D_j\right\Vert_1+1\leq\left\Vert K_n^T\right\Vert_1+4R^3+2.
\end{eqnarray*}
and
	\begin{eqnarray*}
		\Vert V_4\Vert_p
		&\leq&\left\Vert \sum_{j=M_N}^{n}t_{j,n}\left(D_j-D_{M_N}\right)\right\Vert_1 \left\Vert S_{M_{N+1}}\left(f-S_{M_N}f\right)\right\Vert_p\\
		&\leq&\left(\left\Vert K_{n}^{T}\right\Vert_1 +4R^3+2\right) \left\Vert f-S_{M_N}f\right\Vert_p\leq \left(V(\sigma_{n}^T,n)+4R^3+2\right)\omega_p\left(f,1/M_N\right). 
	\end{eqnarray*}
Combining the estimates for $V_1$, $V_2$, $V_3$, and $V_4$ yields the
desired inequality. The proof is complete.
\end{proof}

\section{Further remarks, comments and open problem}

The following result follows immediately from Theorems~\ref{BD4_5} and~\ref{th4}:
\begin{corollary}
Let $G_m$ be a bounded Vilenkin group. Suppose that, for every fixed $k\in\mathbb{P}$, the row
$\{t_{i,k}:1\leq i\leq k\}$ is non-decreasing and that \eqref{28} is satisfied.
\begin{enumerate}
	\item[(a)] For every $f\in L_1(G_m)$,
	$
	\sigma^{T}_{M_k}(f)\to f
	\quad \text{a.e. as }\quad k\to\infty.
	$
	
	\item[(b)] If $\mathbb{E}\subset G_m$ has measure zero, then there exists
	$
	f\in\bigcap_{1\leq p<\infty}L_p(G_m)
	$
	such that the sequence $\{\sigma^{T}_{M_k}(f;x)\}_{k\in\mathbb{P}}$
	diverges at every point $x\in\mathbb{E}$.
\end{enumerate}
\end{corollary}

\begin{proof}
Since almost every point of $G_m$ is a Vilenkin--Lebesgue point of every
$f\in L_1(G_m)$ (see \cite[Corollary~1]{NPTW}), part (a) follows immediately
from Theorem~\ref{BD4_5}.

Part (b) follows from Theorem~\ref{th4} by taking
$
n_k=M_k,\ k\in\mathbb P.
$
\end{proof}

\begin{remark} 
Theorem~\ref{BD4_5} applies to matrix-transform methods which are not covered by the previously known sufficient conditions for convergence. Indeed, consider the triangular matrix defined by $t_{k,n}=0 \ (1\leq k<n), \ t_{n,n}=1.$  Its rows are non-decreasing, condition \eqref{28} is satisfied, and the corresponding matrix-transform method is regular. Indeed, for every fixed $k\in\mathbb P$, we have $t_{k,n}=0$ whenever $n>k$, and hence $t_{k,n}\to0$ as $n\to\infty$. On the other hand, $ t_{n,n}=1,$ and therefore the condition $t_{n,n}=O\left(\frac1n\right),$ used in \cite{BN-MN} for non-decreasing rows, is not satisfied. Moreover, $ \sum_{k=1}^{n}t_{k,n}^{\,2}=1, $ so the condition \eqref{cond-0}, 
used by Blahota and G\'at \cite{BlGa} in the Walsh setting, also fails. Nevertheless, Theorem~\ref{BD4_5} applies and yields \[ \sigma_{M_n}^{T}(f;x)=S_{M_n}(f;x)\to f(x) \] at every Vilenkin--Lebesgue point $x$ of $f$. Thus, Theorem~\ref{BD4_5} provides a complementary convergence result outside the previously known coefficient conditions.
\end{remark}

\begin{remark}
	Theorem~\ref{th4} extends the corresponding divergence results for Fej\'er means \cite{AT2}, N\"orlund means \cite{APST}, and $T$-means \cite{APTY}, the latter two generated by monotone coefficient sequences, to arbitrary regular matrix-transform methods without monotonicity restrictions on the matrix coefficients and along any prescribed increasing sequence of indices.
\end{remark}

\begin{remark}
	Since N\"orlund means generated by non-increasing sequences and
	$T$-means generated by non-decreasing sequences are particular cases
	of matrix-transform means with non-decreasing rows,
	Theorem~\ref{BD4_2} provides a general matrix-transform analogue of
	the approximation estimates \eqref{12} and \eqref{13}.
\end{remark}

\begin{remark}
For $n=M_N$, Theorem~\ref{BD4_3}, together with estimate
\eqref{KT1}, yields an approximation estimate of the same general type
as Theorem~\ref{BD4_2}, although Theorem~\ref{BD4_2} gives a sharper
coefficient for the term involving
$\omega_p\left(f,1/M_N\right)$.
\end{remark}
The following result is an immediate consequence of Theorem~\ref{BD4_3}:

\begin{corollary}
	Let $G_m$ be a bounded Vilenkin group, let $\{n_k\}_{k\in\mathbb N}$ be an increasing sequence of positive integers, and let $f\in L_p(G_m)$ with $1\le p\le\infty.$
	Suppose that, for every fixed $k\in\mathbb P$, the row
	$\{t_{i,k}:1\leq i\leq k\}$ is non-decreasing and \eqref{28} is
	satisfied. If, in addition,
	\[
	\omega_p\left(f,\frac{1}{M_{|n_k|}}\right)
	=
	o\left(
	\frac{1}{V(\sigma^T_{n_k},n_k)}
	\right),
	\qquad k\to\infty,
	\]
	then
	\[
	\|\sigma^T_{n_k}f-f\|_p\to0,
	\qquad k\to\infty.
	\]
\end{corollary}

\begin{proof}
	Set $N_k:=|n_k|$. Since the rows are non-decreasing and normalized,
	for every fixed $j\in\mathbb P$ we have
	$
	t_{j,n_k}\leq {1}/{(n_k-j+1)}\to0,
	\ k\to\infty.
	$
	Hence, for every fixed $l\in\mathbb N$,
	\[
	T_{M_{l+1}-1,n_k}-T_{M_l-1,n_k}
	=
	\sum_{j=M_l}^{M_{l+1}-1}t_{j,n_k}\to0.
	\]
	
	Fix $L\in\mathbb P$. By Theorem~\ref{BD4_3}, split the first sum
	at $L$. The part corresponding to $l<L$ tends to zero as
	$k\to\infty$, while
	\[
	\begin{aligned}
	&\sum_{l=L}^{N_k-1}
	\left(
	T_{M_{l+1}-1,n_k}-T_{M_l-1,n_k}
	\right)
	\omega_p\left(f,\frac1{M_l}\right)
	\\
	&\qquad\leq
	\omega_p\left(f,\frac1{M_L}\right)
	\sum_{l=L}^{N_k-1}
	\left(
	T_{M_{l+1}-1,n_k}-T_{M_l-1,n_k}
	\right)
	\leq
	\omega_p\left(f,\frac1{M_L}\right).
	\end{aligned}
	\]
	Letting first $k\to\infty$ and then $L\to\infty$, we conclude that
	the first term in Theorem~\ref{BD4_3} tends to zero.
	Moreover,
	$
	V(\sigma^T_{n_k},n_k)
	\omega_p\left(f,1/{M_{N_k}}\right)\to0$
	by assumption. 
	Therefore, the conclusion follows from Theorem~\ref{BD4_3}, which completes the proof.
\end{proof}

Theorem~\ref{BD4_3} contains estimates \eqref{MS4} and \eqref{17} as special cases:
\begin{corollary}
	Let $G_m$ be a bounded Vilenkin group, let $f\in L_p(G_m)$ with $1\le p\le\infty$, and let $n\in\mathbb P$. Set $N:=|n|$, so that
	$
	M_N\leq n<M_{N+1}.
	$
	Suppose that, for every fixed $k\in\mathbb P$, the row
	$\{t_{i,k}:1\leq i\leq k\}$ is non-decreasing and \eqref{28}  is
	satisfied. 
	
	a) If, in addition, condition \eqref{MS3} holds, then
	\begin{eqnarray}\label{17c}
		&&\Vert \sigma^T_n f-f\Vert_p \leq(4R^3+1)\sum_{l=0}^{N-1}\left(T_{M_{l+1}-1,n}-T_{M_l-1,n} \right)\omega_p\left(f,1/M_l\right)+O\left(\omega_p\left(f,1/M_{N}\right)\right).
	\end{eqnarray}

b) If, in addition, condition \eqref{cond-1} holds, then
\begin{equation}\label{17a}
\begin{aligned}
\left\|\sigma_n^T(f)-f\right\|_p
\leq{}&
c\sum_{j=0}^{|n|-1}
M_j\left(
\sum_{l=1}^{m_j-1}
t_{(l+1)M_j-1,n}\right)
\omega_p\left(f,\frac{1}{M_j}\right)
+
O\left(
\omega_p\left(f,\frac{1}{M_{|n|}}\right)
\right).
\end{aligned}
\end{equation}
Moreover,
\begin{equation}\label{17b}
\begin{aligned}
\left\|\sigma_n^T(f)-f\right\|_p
\leq{}&
c\sum_{j=0}^{|n|-1}
\frac{M_j}{M_{\vert n\vert}}
\omega_p\left(f,\frac{1}{M_j}\right)
+
O\left(
\omega_p\left(f,\frac{1}{M_{|n|}}\right)
\right).
\end{aligned}
\end{equation}
\end{corollary}
\begin{proof}
Suppose first that $\{t_{i,k}:1\leq i\leq k\}$ is non-decreasing and satisfies conditions \eqref{28} and \eqref{MS3}. Then, by Lemma~\ref{weak00},
$
V(\sigma_n^T,n)=O(1),
$
which yields estimate \eqref{17c}.

On the other hand, suppose that ${t_{i,k}:1\leq i\leq k}$ is non-decreasing and satisfies conditions \eqref{28} and \eqref{cond-1}. Again, by Lemma~\ref{weak00},
$
V(\sigma_n^T,n)=O(1).
$
Moreover, 
$t_{k,n}\leq t_{(l+1)M_j-1,n},$  
$lM_j\leq k\leq (l+1)M_j-1.$
Hence
$
\sum_{k=lM_j}^{(l+1)M_j-1}t_{k,n}
\leq
M_j t_{(l+1)M_j-1,n},$
and consequently
\[
T_{M_{j+1}-1,n}-T_{M_j-1,n}=
\sum_{l=1}^{m_j-1}
\sum_{k=lM_j}^{(l+1)M_j-1}t_{k,n}
\leq
M_j\sum_{l=1}^{m_j-1}
t_{(l+1)M_j-1,n}.
\]
Thus, inequality \eqref{17a} follows.

Using the condition \eqref{cond-1} we also get that
$
t_{(l+1)M_j-1,n}\leq t_{n,n}\leq {C_0}/{n}.
$
Consequently,
\[
\begin{aligned}
M_j\sum_{l=1}^{m_j-1}t_{(l+1)M_j-1,n}\leq
M_j(m_j-1)\frac{C_0}{n}\leq
C_0(R-1)\frac{M_j}{n}\leq
C_0(R-1)\frac{M_j}{M_{|n|}}\leq
C\frac{M_j}{M_{|n|}}.
\end{aligned}
\]
If we insert this estimate in \eqref{17a} we get \eqref{17b},
which completes the proof.
\end{proof}

\begin{remark}
Estimate~\eqref{17b}, obtained as a consequence of Theorem~\ref{BD4_3},
generalizes the approximation estimates~\eqref{aaa}, \eqref{14}, and~\eqref{16}
for Fej\'er, N\"orlund, and $T$-means, respectively.
\end{remark}

\textbf{Open problem.}
Goginava and Nagy \cite{GN10} introduced the notion of matrix-transform variation for the Walsh--Paley system and used it to obtain results concerning Lebesgue constants, approximation, norm convergence, and almost everywhere convergence of matrix-transform means. 

On the other hand, Blahota, Persson and Tephnadze \cite{BPT2018} established two-sided estimates for the Lebesgue constants of partial sums with respect to bounded Vilenkin systems in terms of the two variation-type quantities $v(n)$ and $v^{*}(n)$, and applied these estimates to study boundedness and convergence properties of Vilenkin--Fourier partial sums.

Motivated by these results, it is natural to seek an appropriate notion of matrix-transform variation for bounded Vilenkin systems that reflects the $m$-adic structure, incorporates the $v(n)$- and $v^{*}(n)$-type contributions, and reduces to the variation of Goginava and Nagy \cite{GN10} in the Walsh case. In particular, suitable two-sided estimates for $V(\sigma_n^T,n)$ would strengthen Theorem~\ref{BD4_3} and yield sharper approximation and convergence criteria.

\end{document}